\documentclass[11pt]{article}
\usepackage{amsmath, amsfonts, amssymb}
\usepackage{graphicx, amsthm}
\usepackage{color}
\usepackage{cite}

\usepackage{dsfont}
\newcommand{\R}{{\mathbb R}}

\providecommand{\U}[1]{\protect \rule{.1in}{.1in}}

\providecommand{\U}[1]{\protect \rule{.1in}{.1in}}
\providecommand{\U}[1]{\protect \rule{.1in}{.1in}}

\newtheorem{theorem}{Theorem}

\newtheorem{lemma}[theorem]{Lemma}
\newtheorem{proposition}[theorem]{Proposition}
\theoremstyle{definition}
\newtheorem{definition}[theorem]{Definition}

\begin{document}

\begin{center}
{\LARGE Pathological subgradient dynamics}
\end{center}

\medskip

\begin{center}
{\Large \textsc{Aris Daniilidis, Dmitriy Drusvyatskiy}}

\bigskip
\end{center}

\noindent \textbf{Abstract.} We construct examples of Lipschitz continuous functions, with pathological subgradient dynamics both in continuous and discrete time. In both settings, the iterates generate bounded trajectories, 
and yet fail to detect any (generalized) critical points of the function. 
%\vspace{1in}

\vspace{0.55cm}


\vspace{0.55cm}


\section{Introduction}
The subgradient method plays a central role in large-scale optimization and its numerous applications. %Though the performance guarantees of the method are well-understood for smooth and convex problems, much less is known outside of these two problem classes. 
The primary goal of the method for nonsmooth and nonconvex optimization is to find generalized critical points. For example, for a locally Lipschitz continuous function $f$, we may be interested in finding a point $x$ satisfying the inclusion $0\in \partial  f(x)$, where the symbol $\partial  f(x)$ denotes the Clarke subdifferential.\footnote{The subdifferential $\partial  f(x)$ is the convex hull of all limits of gradients taken at points approaching $x$, and at which the function is differentiable.}  %%Generalized critical points of smooth functions $f$ are, of course, simply the critical points
%in the classical sense. The more general theory is particularly interesting to optimization specialists, because critical %points of convex functions are just minimizers, and more generally, a point is critical exactly when an appropriately  %defined directional derivative is nonnegative in every direction.
The main difficulty in analyzing subgradient-type methods is that it is unclear how to construct a Lyapunov potential for the iterates when the target function  is merely Lipschitz continuous. One popular strategy to circumvent this difficulty is to pass to continuous time where a Lyapunov  function may be more apparent. Indeed, for reasonable function classes, the objective itself decreases along the continuous time subgradient trajectories of the function. For example, this is the path classically followed by Bena{\"\i}m et al. \cite{BHS,BHS2}, Borkar  \cite{Borkar}, Ljung \cite{Ljung}, and more recently by Davis et al. \cite{DDKL} and Duchi-Ruan \cite{duchi_ruan}. \smallskip

Setting the formalism, consider the task of minimizing a Lipschitz continuous function $f$ on $\R^d$ by the subgradient method. It is intuitively clear that the asymptotic performance of the algorithm is dictated by the long term behavior of the absolutely continuous  trajectories $\gamma\colon[0,\infty)\to\R^d$ of the associated subgradient dynamical system  
\begin{equation}\label{eqn:subgrad_sys}
	-\dot{\gamma}(t)\in \partial  f(\gamma(t))\qquad \textrm{for a.e. } t\in [0,\infty).
\end{equation}
%Aside from a technical condition,\footnote{The technical condition is that the set Clarke regular values $f(\R^d\setminus (\partial  f)^{-1}(0))$ is dense in $\R$. This property is classical and we do not focus on it here.} 
Asymptotic convergence guarantees for the subgradient method, such as the seminal work of Nurminskii \cite{Nurminskii1973} and Norkin \cite{norkin1978nonlocal}, and the more recent works of Duchi-Ruan \cite{duchi_ruan} and Davis et al. \cite{DDKL}, rely either explicitly or implicitly on the following assumption.
\begin{enumerate}
%\item (Weak Sard) The set complement of Clarke critical values of $f$ is dense.
\item[-] (Lyapunov) For any trajectory $\gamma(\cdot)$ emanating from a noncritical point of $f$, the composition $f\circ \gamma$ must strictly decrease on some small interval $[0,T)$.
\end{enumerate}
 For example, it is known that the Lyapunov property holds for any convex \cite{bruck,Brezis}, subdifferentially regular \cite{DDKL,MMM}, semi-smooth, and Whitney stratifiable functions \cite{DDKL}. Since this property holds for such a wide class of functions, it is natural to ask the following question. 
\begin{center}
-- Is the Lyapunov property simply true {\em for all} Lipschitz continuous functions\,?
\end{center}
In this work, we show that the answer is negative ({\em c.f.} Subsection~\ref{ss3.1}).
Indeed, we will show that there exist pathological Lipschitz functions that generate subgradient curves \eqref{eqn:subgrad_sys} with surprising behavior. As the first example, we construct (see Proposition~\ref{prop3}) a Lipschitz continuous function $f\colon\R^2\to\R$ and a subgradient curve $\gamma\colon [0,T)\to\R^2$ emanating from a non-critical point $\gamma(0)$, such that $f\circ \gamma$ {\em increases $\alpha$-linearly}: 
	$$f(\gamma(t))-f(\gamma(0))\geq \alpha t\,,\qquad \textrm{for all }t\in [0,T].$$
	In particular, the Lyapunov property clearly fails.	
Our second example presents a Lipschitz function $f\colon\R^4\to\R$ and a {\em periodic} subgradient curve $\gamma\colon [0,\infty)\to\R^d$ that contains no critical points of $f$. In particular, the limiting set of the trajectory $\gamma(\cdot)$ is disjoint from the set of critical points of the function $f$ (see Theorem~\ref{Thm_A} in Subsection~\ref{ss3.2}).
	Our final example returns to the discrete subgradient method with the usual nonsummable  but square summable stepsize. We construct a Lipschitz function $f\colon\R^4\to\R$, for which the subgradient  iterates form a limit cycle that is disjoint from the critical point set of $f$ ({\em c.f.} Subsection~\ref{ss3.3}). Thus the method fails to find any critical points of the constructed function. \smallskip

The examples we construct are built from so-called ``interval splitting sets'' (Definition~\ref{def2}). These are the subsets of the real line, whose restriction to any interval is neither zero- nor full-measure.  Splitting sets have famously been used by Ekeland-Lebourg  \cite{Lebourg} and Rockafellar \cite{rockafellar1981favorable} to construct a pathological Lipschitz function, for which the Clarke subdifferential is the unit interval $\partial  f^{\circ}=[-1,1]$ everywhere. Later, it was shown that functions with such pathologically large subdifferentials are topologically \cite{BW2000}  and algebraically  generic \cite{DF2019} (see Section~\ref{s2}). Notice the function above does not directly furnish a counterexample to the Lyapunov property, since every point in its domain is critical. Nonetheless, in this work, we borrow the general idea of using splitting sets to define Lipschitz functions with large Clarke subdifferentials. The pathological subgradient dynamics then appear by an adequate selection of subgradients that yields a smooth vector field with simple dynamics. It is worthwhile to note that  in contrast to the aforementioned works, the functions we construct trivially satisfy the conclusion of the Morse-Sard theorem: the set of the Clarke critical values has zero measure.

Although this work is purely of theoretical interest, it does identify a limitation of the subgradient method and the differential inclusion approach in nonsmooth and nonconvex optimization. In particular, this work supports the common practice  of focusing on alternative techniques (e.g. smoothing \cite{ermolievnorkinwets95,nestSpok17}, gradient sampling \cite{BLO05}) or explicitely targeting better behaved function classes (e.g. weakly convex \cite{semiconcave,Nurminskii1973,MR520481}, amenable \cite{amen}, prox-regular \cite{prox_reg},  generalized differentiable \cite{norkin1978nonlocal,ermoliev2003solution,MR0461556}, semi-algebraic \cite{BDLS,MR2486055}). 

\section{Notation}
\label{s2}
Throughout, we let $\R^d$ denote the standard $d$-dimensional Euclidean space with inner product $\langle \cdot,\cdot\rangle$ and the induced norm $\|x\|=\sqrt{\langle x,x\rangle}$. The symbol $\mathbb{B}$ will stand for the closed unit ball in $\R^d$.
For any set $A\subset\R^d$, we let $\chi_A$ denote the function that evaluates to one on $A$ and to zero elsewhere. Throughout, we let $m(\cdot)$ denote the Lebesgue measure in $\R$.

A function $f\colon\mathcal{U}\to \mathbb{R}$, defined on an open set $\mathcal{U}\subset\R^d$, is
called {\em Lipschitz continuous} if there exists a real $L>0$ such that the estimate holds:
\begin{equation}\label{eqn:lip_est}
|f(x)-f(y)|\leq L\| x-y\|\qquad \forall x,y\in \mathcal{U}.
\end{equation}
The infimum of all constants $L>0$ satisfying \eqref{eqn:lip_est} is called the Lipschitz modulus of $f$ and will be denoted by $\|f\|_{\mathrm{Lip}}$.
%In this case
%the Lipschitz constant of $f$ is defined as follows:
%\begin{equation}
%||f||_{\mathrm{Lip}}=\sup_{x,y\in \mathcal{U}\,,\,x\neq y}\left \{
%\frac{|f(x)-f(y)|}{\Vert x-y\Vert}\right \}  =\inf \left \{
%k>0\,:\,|f(x)-f(y)|\leq k\Vert x-y\Vert,\, \text{for all }x,y\in
%\mathcal{U}\right \}  \label{LipConst}
%\end{equation}
For any  Lipschitz continuous function $f$, we  let $\mathcal{D}_{f}\subset \mathcal{U}$ denote the set of points where $f$ is differentiable. The classical Rademacher's
theorem guarantees that $\mathcal{D}_{f}$ has full Lebesgue measure in  $\mathcal{U}$. The {\em Clarke subdifferential of $f$ at $x$} is then defined to be the set
\begin{equation}
\partial f(x):={\mathrm{conv}}  \left \{  \lim_{k\to
\infty}\ \nabla f(x_{k}): x_k\to x  \textrm{ and }\{x_{k}\}\subset \mathcal{D}_{f}\right \}, \label{Clarke}
\end{equation}
where the symbol, ${\mathrm{conv}}$, denotes the convex hull operation. It is important to note that in the definition, the set $\mathcal{D}_{f}$ can be replaced by any full-measure subset $\mathcal{D}\subset\mathcal{D}_{f}$;
see \cite[Chapter~2]{Clarke}  for details. It is easily seen that $\partial f(x)$ is a nonempty convex compact set, whose elements are bounded in norm by  
$\|f\|_{\mathrm{Lip}}$. A point $x$ is called {\em (Clarke) critical} if the inclusion $0\in \partial f(x)$ holds. We will denote the set of all critical points of $f$ by $\mathrm{Crit}(f)$. A real number $r\in\R$ is called a {\em critical value} of $f$ whenever there exists some point $x\in \mathrm{Crit}(f)$ satisfying $f(x)=r$.

Though the definition of the Clarke subdifferential is appealingly simple, the behavior of the (set-valued) map $x\rightrightarrows \partial  f(x)$ can be quite pathological;  e.g. \cite{BMW1997,CPT2005,Dymond-Kaluza2019}. For example, there exists a $1$-Lipschitz function $f\colon\R\to\R$ having maximal possible subdifferential $\partial f(x)=[-1,1]$ at every point on the real line. The first example of such {\em Clarke-saturated} functions appears in \cite[Proposition~1.9]{Lebourg}, and is based on interval splitting sets.

\begin{definition}[Splitting set]\label{def2} A measurable set $A\subseteq \mathbb{R}$ is said to {\em split intervals} if for every nonempty interval $I \subset \R$ it holds
\begin{equation}
0<m(A\cap I)<m(I)\,, \label{split}
\end{equation}
where $m$ denotes the Lebesgue measure. 
\end{definition}

%\begin{remark}[Historical commentary]
 The first definition and construction of a splitting
set goes back to \cite{Kirk}, while the first examples of Clarke saturated functions can be found in \cite{rockafellar1981favorable,Lebourg}. The basic construction proceeds as follows. For any fixed set $A\subseteq \mathbb{R}$ that split intervals, define the univariate function
$$f(t)=\int^t_{0} \chi_{A}(\tau)-\chi_{A^c}(\tau)~d\tau\qquad \textrm{for all }t\in \R.$$ 
An easy computation shows the $f$ is Clarke saturated, namely $\partial f(t)=[-1,1]$ for all $t\in\R$. We will use this observation throughout.
Famously, the papers \cite{BMW1997,BW2000} established that, for the uniform topology, a ``generic'' Lipschitz function (in
the Baire sense) is Clarke saturated. Although Clarke-saturation is not
generic under the $\|\cdot\|_{\mathrm{Lip}}$-topology, it has been recently
proved in \cite{DF2019} that the set of Clarke-saturated functions contains a nonseparable Banach space.\footnote{Consequently, in the notation of \cite{G1966}, the set of Clarke-saturated functions is ``spaceable''. Moreover, if we endow the space of Lipschitz continuous functions with the $\|\cdot\|_{\mathrm{Lip}}$-seminorm, the space of pathological functions contains an isometric copy of $\ell^{\infty}$.}
We refer to \cite{AGS2005,BO2014,BPS2014,EGS2014} for recent results on the topic.

%\smallskip
%
%A particular feature of these functions is that every point is Clarke critical
%(i.e. $0\in$ $\partial f(x),$ for all $x\in \mathcal{U}$), therefore,
%the above discussion also exhibits large classes of Lipschitz functions
%failing to satisfy the Morse-Sard theorem. Let us recall that particular
%subclasses of Lipschitz functions, such as semialgebraic/tame Lipschitz
%functions (\cite{BDLS}) or finite continuous selections of smooth functions
%(\cite{BDD}) do satify the conclusion of the Morse-Sard theorem (i.e. the set
%$f(\mathrm{Crit}_{f})$ is Lebesgue null, where $\mathrm{Crit}_{f}
%:=\{x\in \mathcal{U}:0\in \partial f(x)\}$), and the size of their
%subdifferential is generically small (see for instance \cite{DL2013} for the
%semialgebraic case).
%
%\smallskip
%
%In this work we show that (first order) subgradient systems may fail to detect
%critical points in full generality. It is interesting to note that this
%failure appears to functions with large Clarke subdifferentials, yet with a control on their
%generalized critical values (they satisfy Morse-Sard theorem).

\section{Main results}

In this section we construct the three pathological examples announced in
the introduction. Our first example will make use of a splitting set that
satisfies an auxiliary property. The construction is summarized in
Lemma~\ref{L1}. The proof is essentially standard, and therefore we have placed it in the Appendix. 
%\ref{sec:appendix}.

\begin{lemma}[Controlled splitting]
\label{L1} For every $\lambda\in(\frac{1}{2},1)$, there exists a measurable set $
A\subset {\mathbb{R}}$ that splits intervals and satisfies:\newline
\begin{equation}
m(A\cap \lbrack 0,t])\geq \lambda t\,,\quad \text{for every }t>0.  \label{A}
\end{equation}
\end{lemma}

\subsection{Nondecreasing subgradient trajectories}
\label{ss3.1}

The following proposition answers to the negative the first question of the introduction, revealing that the Lyapunov property for a subgradient trajectory may fail.  

\begin{proposition}
[Linear increase along orbits] \label{prop3}Let $\alpha > 0$ be arbitrary. Then, there exists a
Lipschitz continuous function $f\colon\mathbb{R}^{2}\rightarrow \mathbb{R}$ and a
subgradient orbit $\gamma:[0,+\infty)\rightarrow \mathbb{R}^{d}$ emanating from
a noncritical point, meaning
\begin{equation}
\left \{
\begin{aligned}
-\dot{\gamma}(t)&\in \partial f(\gamma
(t))~\textrm{for a.e. }t\in [0,+\infty),\\
\gamma(0)&\notin \mathrm{Crit}(f)
\end{aligned}
\right.  \label{1}
\end{equation}
and satisfying the linear increase guarantee
\[
f(\gamma(t))-f(\gamma(0))\geq \alpha \,t\,,\quad \text{for every }t\in \lbrack0,+\infty).
\]

\end{proposition}

\noindent \textbf{Proof.} According to Lemma~\ref{L1}, there exists a constant
$\lambda>\frac{1}{2}$ and a set $A\subset \R$ satisfying \eqref{A}. Define the constant $\mu:=\sqrt{\frac{\alpha+1}{2\lambda-1}}$ and define the function $f\colon\R^2\to\R$ by
\[
f(x,y):=-x+\mu \int \limits_{0}^{y} \, \left(  \chi_{A}(\tau)-\chi_{A^{c}}
(\tau)\right)  \,d\tau.
\]
It is easily seen that $f$ is Lipschitz continuous and the Clarke subdifferential of $f$ is given by
\[
\partial f(x,y)=\{-1\} \times \lbrack-\mu,\mu],\quad \text{for every
}(x,y)\in \mathbb{R}^{2}.
\]
Notice that  $\partial f(x,y)$ contains the
direction $u=-(1,\mu)$ at every point $(x,y)$. Taking into account $\mathrm{Crit}(f)=\emptyset$, we deduce that the
curve $\gamma(t)=-tu=(t, \mu t)$ satisfies the system (\ref{1}). Moreover, we have the estimate
\[
f(\gamma(t))=f(t,\mu t)=-t+\mu\, \left(  m(A\cap \lbrack0,\mu t])-m(A^{c}\cap
\lbrack0, \mu t])\right)  \geq \left( (2\lambda-1)\mu^{2}-1\right) t\, \geq
\alpha \,t\,.
\]
The proof is complete.$\qquad \hfill \Box$

%%%%%%%%%%%%%%%%%%%%%%%%%%%%%%%%%%%%%%%%%%%%%%%

\subsection{(Periodic) subgradient orbits without critical points}
\label{ss3.2}

Next, we present the second example announced in the introduction, namely a Lipschitz continuous function along with a periodic subgradient curve that contains no critical points of the function. We begin with the following intermediate construction. Henceforth, the symbol $\mathbb{B}_{\infty}$ will denote the closed unit $\ell_{\infty}$-ball in $\R^d$, and  $i$ will denote the imaginary unit.

\begin{theorem}\label{thm:basic_period}
Fix an arbitrary real $M>0$ and $b\in (0,\frac{M}{2})$, and let $A\subset\R$ be a measurable subset that splits intervals. Define the function $\Phi\colon\R^2\to\R$ by
$$\Phi(x,y):=xy+M\, \int \limits_{0}^{y} \left(  \chi_{A}(\tau
)-\chi_{A^{c}}(\tau)\right)  \,d\tau.$$
Then the following are true:
\begin{itemize}
\item[\rm{(i).}] The function $\Phi$ is $2 M$-Lipschitz continuous when restricted to the ball $b\mathbb{B}_{\infty}$.
\item[\rm{(ii).}] Equality holds: $b\mathbb{B}_{\infty}\cap \mathrm{Crit}(\Phi)=[-b,b]\times \{0\}$.
\item[\rm{(iii).}] For any real $r<b$ and $\theta\in\R$, the curve $\gamma(t)=re^{i(t+\theta)}$
is a subgradient orbit of $\Phi$, that is $-\dot{\gamma}(t)\in
\partial \Phi(\gamma(t))$ for all $t> 0$.
\end{itemize}
\end{theorem}
\noindent \textbf{Proof.} The standard sum rule yields the expression for the subdifferential 
\begin{equation}
\partial \Phi(x,y)=\left \{ (y,x+h):\,h\in \lbrack-M,M]\, \right \}
\,,\quad \text{for all }(x,y)\in \R^2. \label{subd-Phi}
\end{equation}
The first claim then follows immediately by noting 
$$\max_{(x,y)\in b\mathbb{B}_{\infty}}\max_{v\in \partial \Phi(x,y)}\|v\|\leq 2M.$$
The second claim also follows immediately from the expression \eqref{subd-Phi}. Next, for any $(x,y)\in b\mathbb{B}_{\infty}$ we can take
$h=-2x$ in (\ref{subd-Phi}),  yielding the selection $(-y,x)\in-\partial \Phi(x,y).$
Thus the inclusion $\dot{\gamma}\in-\partial \Phi(\gamma)$ holds as long as the curve $\gamma(t) =(x(t),y(t))$ satisfies the ODE
\begin{equation}\label{eqn:ode_sat}
\dot{x}(t)=-y(t),~~
\dot{y}(t)=x(t)\qquad \forall t.
\end{equation}
Clearly, the curve $\gamma(t):=re^{i(t+\theta)}$ indeed satisfies $\eqref{eqn:ode_sat}$.
%Trivially the solutions of this ODE are precisely the curves satisfying 
%\eqref{periodic} with $c=\|\gamma(0)\|$. This completes the proof.
\qquad $\hfill \Box$

\bigskip

Thus Theorem~\ref{thm:basic_period} provides an example of a periodic subgradient curve $\gamma(\cdot)$ for a Lipschitz continuous function $\Phi\colon\R^d\to\R$. The deficiency of the construction is that $\gamma$ does pass through some critical points of $f$. We will now see that by doubling the dimension, we can ensure that the periodic curve never passes through the critical point set.

\begin{theorem}\label{thm:thm2}
(Periodic subgradient orbits without critical points)\label{Thm_A} There exists
a Lipschitz continuous function $f\colon\mathcal{U}\longrightarrow \mathbb{R}$, defined on an open set $\mathcal{U}\subset \mathbb{R}^{4}$, and a periodic analytic curve $\gamma\colon \R\to \mathcal{U}$, satisfying 
$$\dot{\gamma}(t)\in-\partial f(\gamma(t))\quad \forall t\qquad \textrm{and}\qquad \operatorname{Im}(\gamma)\, \cap \, \mathrm{Crit}(f)=\emptyset.$$ 
\end{theorem}

\begin{figure}[h!]
	\centering
	\includegraphics[scale=0.8]{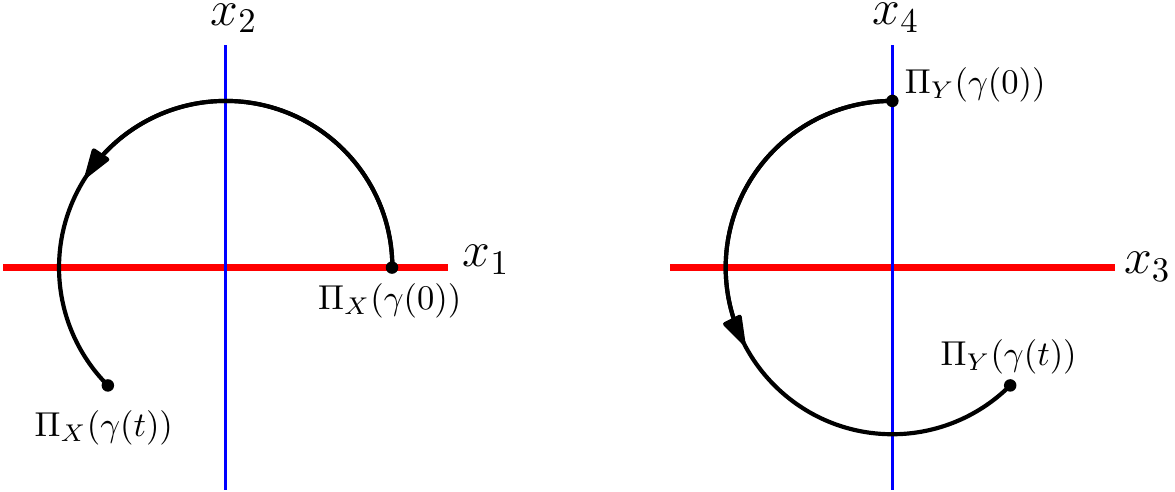}
	\caption{Subgradient curve $\gamma(\cdot)$ generated by $f$ is depicted in black; the set of critical points of $f$ is depicted in red. The picture on the left shows the projection of $\gamma$ onto the coordinate $(x_1,x_2)$, while the picture on the right shows the projection of $\gamma$ onto the coordinate $(x_3,x_4)$.}
	\label{fig:fig2}
\end{figure}

\noindent \textbf{Proof.} Let $M$, $b$, $A$, and $\Phi$ be as defined in Theorem~\ref{thm:basic_period}.
Set $\mathcal{U}:=bB_{\infty}\times bB_{\infty}$ and define the function
\begin{equation}
\left \{
\begin{array}
[c]{l}
f\colon\mathcal{U}\to \mathbb{R}\vspace{0.3cm}\\
f(x_{1},x_{2},x_{3},x_{4})=\Phi(x_{1},x_{2})+\Phi(x_{3},x_{4})
\end{array}.
\right.  \label{f}
\end{equation}
It follows easily from Theorem~\ref{thm:basic_period} that the critical point set is given by
\begin{equation}\label{eqn:crit_set_def}
\mathrm{Crit}(f)=\left \{  (x_{1},0,x_{3},0):\, \max \{|x_{1}|,|x_{3}|\} \leq
b\, \right \}.
\end{equation}
Define the curve $\gamma\colon \R\to \mathcal{U}$ by
$\gamma(t)=\frac{b}{2}\left(e^{it},e^{i\left(t+\frac{\pi}{2}\right)}\right).$
Theorem~\ref{thm:basic_period} immediately guarantees the inclusion $\dot\gamma(t)\in -\partial  f(\gamma(t))$ for all $t\in \R$, while the expression \eqref{eqn:crit_set_def} implies $\operatorname{Im}(\gamma)\, \cap \, \mathrm{Crit}(f)=\emptyset$. See Figure~\ref{fig:fig2} for an illustration.
$\hfill \Box$

\bigskip

It is worthwhile to note that the function $f$ in Theorem~\ref{Thm_A} trivially satisfies the conclusion of the Morse-Sard theorem, since $f(\mathrm{Crit}(f))=\{0\}$.
%
%\bigskip
%
%\noindent \textbf{Remarks. (i). }Let $\beta>0$ be any constant. Then taking
%$M=\frac{2\beta}{\sqrt{10}},$ $b=\frac{\beta}{\sqrt{10}}=\frac{M}{2}$ and
%extending $\Phi$ to $\mathbb{R}^{2}$ by the McShane formula (inf-convolution
%with $\beta \,d(\cdot,\mathbf{\bar{B})}$) we deduce that the function $\Phi$
%can be defined in $\mathbb{R}^{2}$ and $||\Phi||_{\mathrm{Lip}}=\beta$.
%Therefore, the function $f$ of Theorem~\ref{Thm_A} can be taken to be
%Lipschitz continuous of a prescribed constant. \smallskip
%
%
%\textbf{(ii).} The function $f:\mathcal{U}\rightarrow \mathbb{R}$ of (\ref{f})
%satisfies trivially the Morse-Sard theorem. Moreover, for the choice
%$\beta=\frac{\sqrt{2}}{2},$ $M=\frac{2}{\sqrt{10}}\beta=\frac{\sqrt{5}}{5}$
%and $b=\frac{\sqrt{5}}{10}$ we get $||f||_{\mathrm{Lip}}=1.$ Moreover, for
%\[
%\mathcal{G}=\left \{  x\in \mathcal{U}:\,|x_{1}|+|x_{3}|\neq0\; \&
%\;|\frac{x_{2}}{x_{1}}|\neq|\frac{x_{4}}{x_{3}}|\; \right \}
%\]
%we observe that the subgradient orbit (\ref{3}) with initial point
%$\gamma(0)\in \mathcal{G}$ is periodic and free of critical points.
%Consequently, the conclusion of Theorem \ref{Thm_A} can be reinforced as follows:
%\begin{center}
%-- For a generic initial point there exists a period subgradient orbit
%free of critical points.
%\end{center}

\subsection{Subgradient sequences without reaching critical points}
\label{ss3.3}
We next present the final example announced in the introduction. We exhibit a Lipschitz continuous function $f$ such that
the subgradient method, which can access $f$ only by querying sugradients, fails to detect critical points in any sense, under any choice of 
(nonsummable, square summable) steps $\{t_n\}_{n\geq 1}$.
As the initial attempt at the construction, one may try applying the subgradient method to the function $f$ constructed in Theorem~\ref{thm:thm2}, since it has periodic subgradient orbits in continuous time. The difficulty is that when applied to this function, the subgradient iterates (in discrete time) quickly grow unbounded. Therefore, as part of the construction, we will modify the function $f$ from Theorem~\ref{thm:thm2} by exponentially damping its slope.

\begin{theorem}
[Subgradient method does not detect critical points]\label{Thm_B} There exists a Lipschitz continuous function $f\colon\mathbb{R}^{4}\to \mathbb{R}$, a subgradient selection $G(X)\in \partial f(X),$ and initial condition $\bar{X}
\in \mathbb{R}^{4}$ such that
\begin{equation}
\text{for every sequence\ }\{t_{n}\}_{n}\subset(0,+\infty)\; \text{with}\;
\sum_{n\geq1} t_{n}=+\infty \; \text{and} \; \sum_{n\geq1} t_{n}^{2}<+\infty \label{4}
\end{equation}
the subgradient algorithm
\begin{equation}
\left \{
\begin{aligned}
X_{n+1}&=X_{n}-t_{n}G(X_n)\\
X_{1}&= \bar{X}
\end{aligned}
\right\}  \label{5}
\end{equation}
generates a bounded sequence $\{X_{n}\}_{n\geq1}$ whose accumulation points do not meet $\mathrm{Crit}(f)$.
\end{theorem}

\noindent \textbf{Proof.} 
Let us first define the function $\rho:\mathbb{R}
^{2}\to \mathbb{R}$ by
\[
\rho(x,y)=\exp(-y^{2}) \int \limits_{0}^{x} \exp(-\tau^{2})\,d\tau,
\]
and notice that
\[
\nabla \rho(x,y)=\left(  e^{-(x^{2}+y^{2})},\,-2y\cdot\rho(x,y)\right)  .
\]
Then for $\delta>0$ and $M\geq \frac{\delta}{2}(\sqrt{\pi}+1)$, define the function $\phi\colon\mathbb{R}^{2}\to \mathbb{R}$ by
\begin{equation}\label{phi}
\phi(x,y):=\delta \,y\, \rho(x,y)+M\, \int \limits_{0}^{y} \left(  \chi_{A}(\tau
)-\chi_{A^{c}}(\tau)\right)  \,d\tau.
\end{equation}
An easy calculation shows
that $\phi$ is Lipschitz continuous and its subdifferential is given
by
\begin{equation}
\partial \phi(x,y)=\left \{  \delta \,y\,e^{-(x^{2}+y^{2})}\right \}  \,
\times \, \left\{\left (  \delta \,(1-2y^{2})\,e^{-y^{2}}\, \int \limits_{0}^{x}
e^{-\tau^{2}}d\tau \right )  +[-M,M]\right\}. \label{subd-phi}%
\end{equation}
It follows from (\ref{subd-phi}) that a point $(x,y)\in \mathbb{R}^{2}$ is
Clarke critical for $\phi$ if and only if $y=0,$ that is
\[
\mathrm{Crit}(\phi)=\mathbb{R}\times \{0\}.% \quad \text{and\quad}\phi
%(\mathrm{Crit}(\phi))=\{0\}.
\]
%Therefore, $0$ is the only critical value of $\phi$ (i.e. $\phi$ satisfies
%Morse-Sard theorem).
 We claim that for every $(x,y)\in \mathbb{R}^{2}$, we have
\begin{equation}
g(x,y):=\delta \,e^{-(x^{2}+y^{2})}\,(y,-x)\in \partial \phi(x,y).
\label{subg-field}
\end{equation}
To see this, denoting by $\pi_2:\R^2\rightarrow \R$ the projection to the second coordinate, we observe: 
\begin{align*}
\left| \pi_2(g(x,y))- \delta \,(1-2y^{2})\,e^{-y^{2}}\, \int \limits_{0}^{x}
e^{-\tau^{2}}d\tau\right|&=\delta e^{-y^2}\left|(1-2y^{2})\, \int \limits_{0}^{x}
e^{-\tau^{2}}d\tau-xe^{-x^2}\right|\\
&\leq  \delta e^{-y^2}|1-2y^{2}|\cdot\frac{\sqrt{\pi}}{2}+|x|e^{-x^2}\\
&< \delta\cdot\frac{\sqrt{\pi}}{2}+\frac{1}{2}\leq M.
\end{align*}
Notice that the integral curves of the above vector field
\[
(\dot{x},\dot{y})=g(x,y)=\delta \,e^{-(x^{2}+y^{2})}\,(y,-x),
\]
are the homocentric cycles $x^{2}+y^{2}=r^{2}$ for any $r\geq 0.$ %(The reader might want to compare the above with Lemma \ref{Lemma_periodic}.)

$\smallskip$

Consider now applying a subgradient method to $\phi$. Namely, let $\gamma_{n}:=(x_{n},y_{n})$ and consider the subgradient sequence
\begin{equation}
\gamma_{n+1}=\gamma_{n}-t_{n}g_n, \label{6}%
\end{equation}
where we set $g_n:=g(x_n,y_n)\in \partial f(x_n,y_n)$.
 Then the norms
$r_{n}:=\|\gamma_{n}\|$ satisfy
\begin{equation}
\|g_n\|=\delta r_{n}e^{-r_{n}^{2}}\leq \frac{\delta}{\sqrt{2e}}. \label{7}%
\end{equation}
Since $g_n$ is tangent at $\gamma_{n}$ to the homocentric cycle centered at
$0$ with radius $r_{n},$ we deduce easily from (\ref{6}) that the sequence
$\{r_{n}\}_{n\geq1}$ is strictly increasing. On the other hand, by Pythagoras
theorem and (\ref{7}) we deduce:
\[
r_{n+1}^{2}\, = \,r_{n}^{2}\,+\,t_{n}^{2}\,\|g_{n}\|^{2}\leq \,r_{n}
^{2}\,+\, \left(  \frac{\delta^{2}}{2e}\right)  t_{n}^{2},
\]
and by induction
\[
r_{n+1}^{2}\, \leq \,r_{1}^{2}\,+\left(  \frac{\delta^{2}}{2e}\right)
\sum_{n\geq1} t_{n}^{2}<+\infty.
\]
Therefore, $\{r_{n}\}_{n\geq1}$ is bounded and the sequence $\gamma
_{n}:=(x_{n},y_{n})$ has accumulation points.

%\begin{figure}[h!]
%	\centering
%	\includegraphics[scale=0.8]{Fig3.pdf}
%	\caption{Consecutive iterates of the subgradient sequence of $\phi$.}
%	\label{fig:fig3}
%\end{figure}

The proof is not yet complete, since in principle, the accumulation points of the sequence $(x_{n},y_{n})$ may be critical. To eliminate this possibility, we proceed as in the proof of Theorem~\ref{Thm_A} by doubling the dimension.
Namely, define the function
\begin{equation}
\left \{
\begin{array}
[c]{l}
f\colon\mathbb{R}^{4}\to \mathbb{R}\vspace{0.3cm}\\
f(x,y,z,w)=\phi(x,y)+\phi(z,w)
\end{array}
\right. , \label{f-discr}
\end{equation}
and observe that $f$ is Lipschitz continuous and equality holds:
\[
\mathrm{Crit}(f)=\mathbb{R}\times \{0\} \times \mathbb{R}\times \{0\}.
\]

We shall now prescribe a subgradient selection:
\begin{equation}
G(x,y,z,w)=(g(x,y),g(z,w))\in\partial f(x,y,z,w),
\label{G--field}
\end{equation}
where $g(\cdot,\cdot)$ is defined in \eqref{subg-field}. Let us also prescribe the
initial condition
\[
X_{1}\equiv \bar{X}=(1,0,0,1)\in \mathbb{R}^{4}.
\]
Then (\ref{5}) generates a bounded sequence
\[
X_{n}:=(x_{n},y_{n},z_{n},w_{n})\text{\quad with\quad}G_{n}:=G(X_{n}
)\in \partial f(X_{n})
\]
which splits in $\R^2 \times \R^2$ as follows:
\[
X_{n}=(\gamma_{n},\tilde{\gamma}_{n}),\;G_{n}=(g_{n},\tilde{g}_{n})\in
\partial f(X_{n})\text{\quad with\quad}\gamma_{n}:=(x_{n},y_{n}
)\quad \text{and\quad \ }\tilde{\gamma}_{n}:=(z_{n},w_{n})
\]
such that
\begin{equation}
\left \{
\begin{array}
[c]{c}
\gamma_{n+1}=\gamma_{n}-t_{n}g_{n}\\
g_{n}\in \partial \phi(\gamma_{n})
\end{array}
\right\}  \quad \text{and}\quad \text{ }\left \{
\begin{array}
[c]{c}
\tilde{\gamma}_{n+1}=\tilde{\gamma}_{n}-t_{n}\tilde{g}_{n}\\
\tilde{g}_{n}\in \partial \phi(\tilde{\gamma}_{n})
\end{array}
\right\},  \label{rotational}
\end{equation}
respectively. Thanks to the initial condition, for every $n\geq1$ the vector
$\tilde{\gamma}_{n}$ is a $\frac{\pi}{2}$-rotation of the vector
$\gamma_{n}.$ Taking into account this rotational symmetry we deduce
that all limit points of $\{X_{n}
\}_{n\geq1}$ lie outside of $\mathrm{Crit}(f)$. $\hfill \Box$

\bigskip

It is worthwhile to note again that the function $f$ in Theorem~\ref{Thm_B} trivially satisfies the conclusion of the Morse-Sard theorem, since $f(\mathrm{Crit}(f))=\{0\}$.

\vspace{0.8cm}

\textbf{Acknowledgments.} A major part of this work was done during a
research visit of the first author to the University of Washington (July~2019). This author thanks the host institute for hospitality. The first author's research
has been supported by the grants CMM-AFB170001, FONDECYT 1171854 (Chile) and
PGC2018-097960-B-C22, MICINN (Spain) and ERDF (EU). The research of the second author has been supported by the NSF DMS 1651851 and CCF 1740551 awards.

\bigskip
\bibliographystyle{plain}

\vspace{0.8cm}

\noindent Aris Daniilidis 

\medskip

\noindent DIM--CMM, UMI CNRS 2807\newline Beauchef 851, FCFM, Universidad de
Chile \smallskip

\noindent E-mail: \texttt{arisd@dim.uchile.cl}; URL: \texttt{http://www.dim.uchile.cl/{\raise.17ex\hbox{$\scriptstyle\sim$}}arisd}

\medskip

\noindent Research supported by the grants: \newline CMM AFB170001,
FONDECYT 1171854 (Chile), PGC2018-097960-B-C22 (Spain and EU).

\vspace{0.5cm}

\noindent Dmitriy Drusvyatskiy

\medskip

\noindent Department of Mathematics, University of Washington, Seattle, WA 98195 \smallskip

\noindent E-mail: \texttt{ddrusv@uw.edu}; URL:
\texttt{www.math.washington.edu/{\raise.17ex\hbox{$\scriptstyle\sim$}}ddrusv}

\medskip

\noindent Research supported by the NSF DMS 1651851 and CCF 1740551 awards

\bigskip
\begin{center}
\rule{16cm}{1.4pt}

\newpage
%\vspace{1cm}
{\Large {\bf Appendix: Proof of Lemma~\ref{L1}}}
\end{center}
\smallskip

\noindent We will first need the following lemma, outlining a standard construction of a fat Cantor set. We
will impose an extra property, which will play a key role in establishing Lemma~\ref{L1}.

\begin{lemma}[$\protect\alpha $-fat Cantor set with attributes]
\label{L0}Fix any $\alpha \in (0,1)$. Then for every interval $
I=[a,b]\subset \mathbb{R}$ there exists an $\alpha $-fat Cantor subset $
\mathcal{F}_{\alpha }(I)\subset I$, that is, a Cantor-type set of total
measure 
\begin{equation*}
m(\mathcal{F}_{\alpha }(I))=\alpha \,m(I).
\end{equation*}
Moreover, for any $\frac{1}{2}<\lambda <1,$ taking $\alpha >\frac{3\lambda }{
2+\lambda }$ we ensure that 
\begin{equation}
m(\mathcal{F}_{\alpha }(I)\cap \lbrack a,x])\geq \lambda \cdot
m([a,x]),\quad \text{for all }x\in \lbrack a,b].  \label{lary}
\end{equation}
\end{lemma}

\noindent \textbf{Proof.} The construction is standard and is sketched for
the reader's convenience. Fix any $\alpha \in (0,1)$. Let us denote $I^{0}=(a,b),$ $\ell _{0}=m(I^{0})=b-a$
and let us set:
\begin{equation}
\delta _{0}:=\left( \frac{1-\alpha }{2}\right) \cdot \ell _{0}>0.
\label{delta}
\end{equation}
We shall construct the fat-Cantor set by removing, successively, countably many intervals, indexed by a dyadic tree. To this end, we start by removing from our initial set $I^{0}$ the interval 
\begin{equation*}
\Delta ^{0}:=\left[\frac{b+a-\delta _{0}}{2},\frac{b+a+\delta _{0}}{2}\right],
\end{equation*}
that is, an interval of length $\delta _{0}$ centered at $\frac{a+b}{2}$
(the midpoint of $I^{0}$). Then from each of the two remaining intervals $
I_{0}^{1}:=(a,\frac{b+a-\delta _{0}}{2})$ and $I_{1}^{1}:=(\frac{b+a+\delta
_{0}}{2},b)$ we subtract intervals $\Delta _{0}^{1}$ and $\Delta _{1}^{1}$
of length 
\begin{equation}
\delta _{1}:=\frac{1}{2}\cdot \frac{\delta _{0}}{2}  \label{delta1}
\end{equation}
centered at the midpoints of $I_{0}^{1}$ and $I_{1}^{1}$ respectively.
Setting 
\begin{equation*}
m(I_{0}^{1})=m(I_{1}^{1}):=\ell _{1}=\frac{\ell _{0}-\delta _{0}}{2}
\end{equation*}
we observe that 
\begin{equation*}
\frac{m(\Delta ^{0})}{m(I^{0})}=\frac{\delta _{0}}{\ell _{0}}=\frac{1-\alpha 
}{2}\quad \text{and\quad }\frac{m(\Delta _{i}^{1})}{m(I_{i}^{1})}=\frac{
\delta _{1}}{\ell _{1}}=\frac{(\delta _{0}/4)}{(\ell _{0}-\delta _{0})/2}=
\frac{1}{2}\left( \frac{1-\alpha }{1+\alpha }\right) ,\quad \text{for }i\in
\{0,1\},
\end{equation*}
and consequently
\begin{equation}
\frac{m(\Delta _{i}^{1})}{m(I_{i}^{1})}
=\left(\frac{1}{1+\alpha}\right)\frac{m(\Delta ^{0})}{m(I^{0})} \; <\,
\frac{m(\Delta ^{0})}{m(I^{0})} ,\quad \text{for }i\in \{0,1\}.
\label{clef}
\end{equation}
The above relation reveals that in the second step, we subtract a
proportionally smaller part of each of the intervals  $I_{i}^{1}$, $i\in
\{0,1\},$ compared with what we subtract from $I^{0}$ in the first
step.\smallskip \newline
%After the second step, we obtain four intervals $I_{j}^{2},$ $j\in
%\{0,1,2,3\},$ of length 
%\begin{equation*}
%m(I_{j}^{2}):=\ell _{2}=\frac{\ell _{1}-\delta _{1}}{2},\quad \text{for }
%j\in \{0,1,2,3\}.
%\end{equation*}
%Then for $j\in \{0,1,2,3\},$ we subtract from the interval $I_{j}^{2},$ an
%interval $\Delta _{j}^{2}$ of length 
%\begin{equation*}
%\delta _{2}=\frac{1}{2}\cdot \frac{\delta _{1}}{2}=\frac{1}{2^{2}}\cdot 
%\frac{\delta _{0}}{2^{2}}
%\end{equation*}
%centered at the midpoint of $I_{j}^{2}$. Notice that for all $i\in \{0,1\}$
%and $j\in \{0,1,2,3\}$ we have: 
%\begin{equation*}
%\frac{m(\Delta _{j}^{2})}{m(I_{j}^{2})}=\frac{\delta _{2}}{\ell _{2}}=\frac{
%(\delta _{1}/4)}{(\ell _{1}-\delta _{1})/2}<\frac{\delta _{1}}{\ell _{1}}=
%\frac{m(\Delta _{i}^{1})}{m(I_{i}^{1})}.
%\end{equation*}
We continue by induction, subtracting at the step $n$
intervals $\Delta _{k}^{n},$ $k\in \{0,1,\ldots ,2^{n}-1\},$ of length 
\begin{equation*}
\delta _{n}=\frac{1}{2}\cdot \frac{\delta _{n-1}}{2}=\frac{1}{2^{n}}
\cdot \frac{\delta _{0}}{2^{2}}
\end{equation*}
centered at the midpoints of the intervals $I_{k}^{n}$, where 
$ m(I_{k}^{n})=\ell _{n}=(\ell _{n-1}-\delta _{n-1})/2$,
so that for all $k^{\prime }\in \{0,1,\ldots ,2^{n-1}-1\}$ and $k\in
\{0,1,\ldots ,2^{n}-1\}$ we have:
\begin{equation}
\frac{m(\Delta _{k}^{n})}{m(I_{k}^{n})}=\frac{\delta _{n}}{\ell _{n}}=\frac{
(\delta _{n-1}/4)}{(\ell _{n-1}-\delta _{n-1})/2}<\frac{\delta _{n-1}}{\ell
_{n-1}}=\frac{m(\Delta _{k^{\prime }}^{n})}{m(I_{k}^{n})}.  \label{cle}
\end{equation}
The above relation says that at the step $n\geq 2$ we remove a proportionally
smaller part of each of the intervals $I_{k}^{n}$, $k\in \{0,1,\ldots
,2^{n}-1\},$ compared to what we did in the previous step to the intervals 
$I_{k^{\prime }}^{n}$, $k^{\prime }\in \{0,1,\ldots ,2^{n-1}-1\}$. \smallskip 

Let now $\mathcal{F}_{\alpha }(I)$ be the complement of the union of all
extracted intervals. Clearly $\mathcal{F}_{\alpha }(I)$ cannot contain any
interval (\emph{i.e.} it is a Cantor-type set), and in view of \eqref{delta}
its total length is 
\begin{equation*}
m(\mathcal{F}_{\alpha }(I))=m(I)-\left( \delta +2(\frac{1}{2}\cdot \frac{
\delta }{2})+2^{2}(\frac{1}{2^{2}}\cdot \frac{\delta }{2^{2}})+\ldots
\right) =m(I)-\delta \sum_{n\geq 0}\frac{1}{2^{n}}=\alpha \,m(I).
\end{equation*}

\noindent It remains to prove that \eqref{lary} holds. To this end, we start by treating the case 
$x\in \lbrack a,b]\diagdown \mathcal{F}_{\alpha }(I).$ \\
Let us first assume 
\begin{equation*}
x\in \Delta ^{0}:=\left[\frac{b+a-\delta _{0}}{2},\frac{b+a+\delta _{0}}{2}\right].
\end{equation*}
In this case (which is the less favorable case) we have:
\begin{equation}
m([a,x])\leq \frac{b-a+\delta _{0}}{2}=\frac{\ell _{0}+\delta _{0}}{2}\quad 
\text{and\quad }m(\mathcal{F}_{\alpha }(I)\cap \lbrack a,x])=\frac{m(
\mathcal{F}_{\alpha }(I))}{2}=\frac{\alpha \ell _{0}}{2}.  \label{k1}
\end{equation}
Therefore \eqref{lary} holds for $\alpha >3\lambda (2+\lambda )^{-1}$.

Assume now that $x\in \Delta _{0}^{1}\cup \Delta _{1}^{1}.$ By
construction, using \eqref{delta}, \eqref{delta1} and \eqref{clef}, we
deduce that for $i\in \{0,1\},$ the set $\mathcal{F}_{\alpha }(I^{0})\cap
I_{i}^{1}$ is an $\alpha _{1}$-fat Cantor set in $I_{i}^{1}$, denoted $
\mathcal{F}_{\alpha _{1}}(I_{i}^{1}),$ where $\alpha _{1}$ satisfies: 
\begin{equation*}
\frac{1-\alpha _{1}}{2}:=\frac{\delta _{1}}{\ell _{1}}=\left( \frac{1}{
1+\alpha }\right) \left( \frac{\delta _{0}}{\ell _{0}}\right) =\left( \frac{1
}{1+\alpha }\right) \left( \frac{1-\alpha }{2}\right) \quad \text{yielding
that }\,\, \alpha _{1}=\frac{2\alpha }{1+\alpha }>\alpha .
\end{equation*}
The above guarantees that for 
$ x\in \Delta _{0}^{1}\subset I_{0}^{1}=\left(a,\frac{b+a-\delta _{0}}{2}\right)$
we have 
\begin{equation*}
m(\mathcal{F}_{\alpha }(I^{0})\cap \lbrack a,x])=m(\mathcal{F}_{\alpha
_{1}}(I_{i}^{1})\cap \lbrack a,x])\geq \lambda \cdot m([a,x])
\end{equation*}
by the previous step. If now 
\begin{equation*}
x\in \Delta _{1}^{1}\subset I_{1}^{1}=\left(\frac{b+a+\delta _{0}}{2},b\right)
\end{equation*}
then by \eqref{k1} and the fact that $\mathcal{F}_{\alpha _{1}}(I_{1}^{1}):=
\mathcal{F}_{\alpha }(I^{0})\cap I_{1}^{1}$ is an $\alpha _{1}$-fat Cantor
set in $I_{1}^{1}$ with $\alpha _{1}>\alpha $ we deduce:
\begin{align*}
m(\mathcal{F}_{\alpha }(I^{0})\cap \lbrack a,x])&=\frac{m(\mathcal{F}_{\alpha
}(I))}{2}+m\left( \mathcal{F}_{\alpha _{1}}(I_{1}^{1})\cap \left[ \frac{
b+a+\delta _{0}}{2},x\right]\right)\\ &\geq \lambda \cdot \left( m\left( \left[a,\frac{
b+a+\delta _{0}}{2}\right]\right) +m\left( \left[\frac{b+a+\delta _{0}}{2},x\right]\right)
\right) ,
\end{align*}
that is \eqref{lary} holds. Continuing, we deduce that \eqref{lary}
holds for all $x\in \lbrack a,b]\diagdown \mathcal{F}_{\alpha }(I)$.

Let now $x\in \mathcal{F}_{\alpha }(I)$. Then there exists $\{x_{n}\}\subset
\lbrack a,b]\diagdown \mathcal{F}_{\alpha }(I)$ with $\{x_{n}\}\rightarrow x$. Then by \eqref{lary} we have 
\begin{equation*}
m(\mathcal{F}_{\alpha }(I^{0})\cap \lbrack a,x_{n}])\geq \lambda \cdot
m([a,x_{n}])\quad \text{for all }n\geq 1.
\end{equation*}
The result follows by passing to the limit as $n \to \infty$. The proof is complete. $\hfill \Box $

\bigskip

\noindent We are now ready to complete the proof of  Lemma~\ref{L1}.

\bigskip

\noindent \textbf{Proof of Lemma~\ref{L1}.} Let us fix $\lambda\in(\frac{1}{2},1)$ and choose $\alpha ,\theta <1$ (close to $1$) such that 
$$\alpha \cdot \theta >\lambda \qquad \text{and} \qquad \alpha >3\lambda (2+\lambda )^{-1}.$$ For any interval $I=[a,b]$, we set 
\begin{equation*}
I^{+}:=[a,a+\theta (b-a)]\quad \text{and\quad }I^{-}:=[a+\theta (b-a),b]
\end{equation*}
and we define the operators:

\begin{itemize}
\item $\mathcal{T}(I):=\mathcal{F}_{\alpha }(I^{+})$ (partial $\alpha $-fat
Cantor subset of $I$ of measure $\alpha \cdot \theta \cdot m(I)$), and

\item $\mathcal{N}(I)=\mathcal{F}_{\frac{1}{2}}(I^{-})$ (partial $\frac{1}{2}
$-fat Cantor subset of $I$ of measure $\frac{1}{2}(1-\theta )m(I)$).
\end{itemize}

\smallskip 

Notice that it is sufficient to construct $A_{0}\subset \lbrack 0,1]$
(splitting the family of intervals of $(0,1)$) satisfying \eqref{A} for $
x\in \lbrack 0,1)$. Indeed, translating the construction by $n$ we obtain $
A_{n}\subset \lbrack n,n+1)$ and we observe that the set $A=\cup _{n}A_{n}$
has the desired property. \smallskip 

To this end, set $I_{0}=[0,1]$ and consider an enumeration $\{I_{n}\}_{n\geq
1}$ of all strict subintervals $I_{n}:=[p_{n},q_{n}]$ of $[0,1]$ with
rational endpoints $p_{n},$ $q_{n}\in \mathbb{Q}$. Set further $T_{0}:=
\mathcal{T}(I_{0})$ and $N_{0}:=\mathcal{N}(I_{0})$ and notice that $\left(
T_{0}\cup N_{0}\right) $ contains no intervals. Let $I_{1}:=(p_{1},q_{1})$
be the first interval of the above enumeration. Then there exists a
subinterval $I_{1}^{\prime }=[p_{1}^{\prime },q_{1}^{\prime }]$ of $I_{1}$
contained in $I_{0}^{{}}\diagdown \left( T_{0}\cup N_{0}\right) $. We set $
T_{1}=\mathcal{T}(I_{1}^{\prime })$ and $N_{1}=\mathcal{N}(I_{1}^{\prime })$
. Similarly, we find a subinterval $I_{2}^{\prime }\subset I_{2}$ contained
in 
\begin{equation*}
I_{0}^{{}}\diagdown \left[ \left( T_{0}\cup T_{1}\right) \cup \left(
N_{0}\cup N_{1}\right) \right] 
\end{equation*}
and set $T_{2}=\mathcal{T}(I_{2}^{\prime })$ and $N_{2}=\mathcal{N}
(I_{2}^{\prime })$ and continue by induction. Then set 
\begin{equation*}
N=\bigcup_{i=0}^{\infty }N_{i}\quad \text{and}\quad A_{0}=I_{0}\setminus N.
\end{equation*}
It is easily seen that $A_{0}$ splits the family of intervals of $I_{0}=[0,1].$ Let us show that \eqref{A} holds.\newline
Let $t\in \lbrack 0,1]=I_{0}.$ If $t\in I_{0}^{+}=[0,\theta ),$ then since $T_{0}=\mathcal{F}_{\alpha }(I_{0}^{+})$  we conclude by Lemma~\ref{L0} that 
\begin{equation*}
m(A\cap \lbrack 0,t])\geq m(T_{0}\cap \lbrack 0,t])\geq \lambda t.
\end{equation*}
If $t\in I_{0}^{-}=[\theta ,1],$ then 
\begin{equation*}
m(A\cap \lbrack 0,t])\geq m(T_{0}\cap \lbrack 0,t])=m(T_{0})=\alpha \cdot
\theta >\lambda \geq \lambda t.
\end{equation*}
Therefore, $A_{0}$ satisfies \eqref{A} for all $t\in [0,1]$ and consequently, so does the set 
$A:=\bigcup_{n=-\infty }^{+\infty }(A_{0}+n)$ for all $t>0$. $\hfill \Box $
\bigskip

\begin{center}
%\noindent
\rule{8cm}{1.4pt}
\end{center}
\end{document}